\input amstex
\documentstyle{amsppt} 
\TagsOnRight
\magnification = \magstep1
\input prepictex
\input pictex
\input postpictex 
\hsize=34 true pc
\vsize=53 true pc
\voffset=0.1 true in
\hoffset=1.0 true cm  
\overfullrule0pt
\NoBlackBoxes 
\def \n {\noindent}

\def \m {\medskip}
\def \b {\bigskip}

\def\einv{\buildrel E^{-1} \over \rightarrow}
\def\emap{\buildrel E \over \rightarrow}   
\def \End {{\text {End}}}
\def\eqdef{\buildrel\text {def} \over =} 
\def\fsl{\frak s\frak l}
\def\fgl{\frak g\frak l}  
  
\def \g {{\frak g}} 
\def \h {{\frak h}} 

\def \id {{\operatorname{id}}}  
\def \la {\lambda} 

\def \ot {\otimes}

\def \sgn{{\text{sgn}}} 
\def \tr {{\text {tr}}}
\def \u {\underline}

\def \Cent {{\Cal C}} 
\def \nn {n}
\def \rr {r}

\topmatter
\title Derangements and Tensor Powers of Adjoint Modules for $\fsl_\nn$
\endtitle
\author Georgia Benkart
\footnote{Supported in part by NSF Grant \#9970119 \newline
2000 Mathematical Subject Classification  Primary 17B10, Secondary 
05E10  \hfill \ }  \\
Stephen Doty \\
{August 2, 2004} \\
\leftheadtext{GEORGIA BENKART, STEPHEN DOTY}
\rightheadtext {Derangements and Tensor Powers of $\fsl_\nn$}
\endauthor 
\abstract We obtain the decomposition of the tensor space 
$\fsl_\nn^{\otimes k}$ as a module for $\fsl_\nn$, find an explicit 
formula for the multiplicities of its irreducible summands, and (when 
$\nn \ge 2k$) describe the centralizer algebra 
$\Cent=\End_{\fsl_\nn}(\fsl_\nn^{\otimes k})$ and its representations.  
The multiplicities of the irreducible summands are derangement numbers 
in several important instances, and the dimension of $\Cent$ is given 
by the number of derangements of a set of $2k$ elements.  \endabstract 
\endtopmatter \document \b \head Introduction \endhead

\m Weyl's celebrated theorem on complete reducibility says that a
\allowbreak finite-dimen- sional module $X$ for a finite-dimensional
simple complex Lie algebra $\g$ is a direct sum of irreducible
$\g$-modules.  However, to determine an explicit expression for the
multiplicities of the irreducible $\g$-summands of $X$ often is a very
challenging task.  In this note we assume $\g = \fsl_\nn$, the simple
Lie algebra of $\nn \times \nn$ matrices of trace 0 over $\Bbb C$, and
view $\fsl_\nn$ as a $\g$-module under the adjoint action $x \cdot y =
[x,y]$.  We take $X$ to be the $k$-fold tensor power of $\fsl_\nn$.
Using combinatorial methods and results developed in [BCHLLS], we
establish an explicit description of the 
irreducible $\g$-summands of $\fsl_\nn^{\ot k}$ (Theorem 1.17)
and determine an expression for their multiplicities (Theorem 2.2).  As a 
consequence of our formula, we obtain the following results, expressed 
in terms of the number $D_k$ of derangements of $\{1,\dots,k\}$: \ \ 
For $\nn \geq 2k$, the dimension of the space of $\g$-invariants in 
$\fsl_\nn^{\ot k}$ is $D_k$; \ the multiplicity of $\fsl_\nn$ in 
$\fsl_\nn^{\ot k}$ is $D_{k+1}$; and the dimension of the centralizer 
algebra $\Cent = \End_{\g}(\fsl_\nn^{\ot k})$ is $D_{2k}$.

\m In Section 3, we identify the centralizer algebra $\Cent$ with a certain
subalgebra of the walled Brauer algebra $B_{k,k}(\nn)$. This subalgebra
has a basis indexed by derangements of $\{1,\dots, 2k\}$.  We then
give a description (for $\nn\ge 2k$) of the irreducible modules for
$\Cent$, and obtain the ``double centralizer'' decomposition of the tensor
space $\fsl_\nn^{\ot k}$ as a bimodule for $\Cent \times \g$.

\m {\bf Acknowledgment.} \ This is a revised version of a paper
by the same title that appeared in {\it  Journal of Algebraic Combinatorics}
{\bf 16} (2002), 31-42.   In particular 1.15-1.18 in that paper have been
revised in 1.15-1.18 here and a few other related minor changes have been made
in the first line of Section 2 and in Section 3.  
We are grateful to Alberto Elduque for alerting us to 
the mistake in the previous version.

\bigskip
\m
\head \S 1.  The Tensor Product Realization \endhead

\m The general linear Lie algebra $\fgl_\nn = \fsl_\nn \oplus \Bbb C I$ of
all $\nn \times \nn$ complex matrices acts on $\fsl_\nn$ via the adjoint
action, and the identity matrix $I$ acts trivially.  Hence,
there is no harm in assuming that $\g$ is $\fgl_\nn$ rather than
$\fsl_\nn$ acting on $\fsl_\nn^{\ot k}$ in what follows; the results are
exactly the same.  This enables us to label the irreducible summands
by pairs of partitions and to apply known results on the decomposition
of tensor products for $\fgl_\nn$.

\b Let $\h$ denote the Cartan subalgebra of $\g = \fgl_\nn$ of diagonal
matrices, and let $\epsilon_i\,: \h \rightarrow \Bbb C$ be the
projection of a diagonal matrix onto its $(i,i)$-entry.  The
irreducible finite-dimensional $\g$-modules are labeled by their
highest weight, which is an integral linear combination
$\sum_{i=1}^\nn \kappa_i \epsilon_{i}$ with $\kappa_1 \geq \kappa_2 \geq
\cdots \geq \kappa_\nn$.  By letting $\lambda = (\lambda_1 \geq
\lambda_2 \geq \cdots)$ denote the sequence of positive $\kappa_i$ and
$\mu = (\mu_1 \geq \mu_2 \geq \cdots)$ be the partition determined by
the negative $\kappa_i$, we may associate to each highest weight a
pair of partitions $(\lambda, \mu)$.  For example, for $\g =
\fgl_{12}$ the highest weight

$$3 \epsilon_1 + 2 \epsilon_2 + 2 \epsilon_3 + 2\epsilon_4 + \epsilon_5
- 4 \epsilon_{10} - 5 \epsilon_{11} - 5 \epsilon_{12}$$

\n is identified with the pair of partitions $\lambda = (3,2,2,2,1) \vdash 10$
and $\mu = (5,5,4) \vdash 14$.  Therefore, 
the set of highest weights for $\g$-modules is in bijection with the set of pairs of
partitions such that the total number of nonzero parts does not  
exceed $\nn$.

\b Let $V = \Bbb C^\nn$ be the natural representation of $\g = \fgl_\nn$
on $\nn \times 1$ matrices by matrix multiplication.  The dual module
$V^*$ may be identified with $1 \times \nn$ matrices, where the
$\g$-action is by right multiplication by the negative of an element
$x \in \g$.  The matrix product

$$V \ot V^{*} \rightarrow \fgl_\nn = \fsl_\nn \oplus \Bbb C I, \qquad u
\ot w^* \mapsto uw^* \tag 1.1$$

\n is a $\g$-module isomorphism which allows us to identify $\fgl_\nn$
with $V \ot V^{*}$.

\b Let $\{v_1, \dots, v_\nn\}$ denote the standard basis of $V$, where
$v_i$ is the matrix having 1 in the $i$th row and $0$ everywhere else.
Assume $\{v_1^{*}, \dots, v_\nn^{*}\}$ is the dual basis in $V^{*}$, so
that $v_{i}^{*}$ has $1$ in its $i$th column and $0$ elsewhere.  The
{\it contraction mapping} $c\,: \, V \ot V^{*} \rightarrow V \ot
V^{*}$ is defined using the trace by

$$c(u \ot w^{*}) = \tr(uw^{*})\sum_{\ell=1}^{\nn} v_{\ell} \ot
v_{\ell}^{*}. \tag 1.2$$

\n Under the isomorphism in (1.1), $v_\ell \ot v_{\ell}^{*}$ is mapped
to the matrix unit $E_{\ell,\ell} \in \fgl_\nn$.  Therefore, we may
identify the image of $c$ with $\Bbb C I$, and the kernel of $c$ with
$\fsl_\nn$.

\m As $c^{2} = \nn c$, the mapping $p = (1/\nn) c$ is an idempotent.  It
is the projection onto the trivial summand $\Bbb C I$, and $\id-p$ is
the projection onto $\fsl_\nn$.  These idempotents are orthogonal,

$$p(\id-p) = 0 = (\id -p)p,$$

\n and satisfy $\id = (\id-p) + p$. (Here $\id$ is the identity map on
$V \ot V^*$.)

\m In order to identify $\fsl_\nn^{\ot k}$ with a summand of

$$M = V^{\ot k} \ot (V^*)^{\ot k} \cong (V \ot V^*)^{\ot k} \cong
\fgl_\nn^{\ot k}, \tag 1.3$$

\n we define the contraction map $c_{i,j}$ to be the contraction $c$
applied to the $i$th factor of $V^{\ot k}$ and the $j$th factor of
$(V^*)^{\ot k}$ according to
 
$$\multline c_{i,j}(u_1 \ot \cdots \ot u_k \ot w_1^* \ot \cdots \ot
w_k^*) \\ = \tr(u_iw_j^*) \sum_{\ell=1}^\nn u_1 \ot \cdots v_\ell \ot
\cdots \ot v_k \ot w_1^* \ot \cdots \ot v_\ell^* \cdots \ot w_k^*,
\endmultline$$

\n where $v_\ell$ is placed in the $i$th slot of $V^{\ot k}$ and 
$v_\ell^*$ in the $j$th slot of $(V^{*})^{\ot k}$.
As before,  $c_{i,j}^2 = \nn c_{i,j}$, so that
 
$$p_i = \frac{1}{\nn} c_{i,i} \tag 1.4$$

\noindent is an idempotent.   

\b
\proclaim {Proposition 1.5} $\ker p_1 \cap \ker p_2 \cap \cdots 
\cap \ker p_k = (\id-p_1)(\id-p_2) \cdots (\id-p_k)M$.  \endproclaim 
\b  
\demo {Proof} The idempotents $p_i$ commute and satisfy
$p_i(\id-p_i) = 0$.    For $J$ a subset of $\{1, \dots, k\}$,
let $p_J = \prod_{j \in J} p_j$.   Set $q_j = \id-p_j$ 
and $q_J = \prod_{j \in J} q_j$.
Then 

$$M = \bigoplus_{J \subseteq \{1, \dots, k\}} p_{J^c} q_J M,$$ 

\n where $J^c = \{1,\ldots,k\}\setminus J$.  This can be argued by
induction on $k$.  Note that the sum is direct because for any fixed
choice of subset $J'$, the idempotent $p_{{J'}^c}q_{J'}$ acts as the
identity on $p_{{J'}^c}q_{J'}M$ and annihilates the remaining terms
$p_{J^c} q_J M$ with $J \neq J'$.  Whenever $j \in J^c$, then $p_{J^c}
q_J M$ is not contained in $\ker p_j$.  Therefore, from the
decomposition of $M$ above, it is easy to see that $\ker p_1 \cap \ker
p_2 \cap \cdots \cap \ker p_k = (\id-p_1)(\id-p_2) \cdots (\id-p_k)M$.
\qed \enddemo

\b 
Henceforth, let 

$$e = (\id-p_1)(\id-p_2) \cdots (\id-p_k)\tag 1.6$$ 

\n so that 

$$eM \cong \fsl_\nn^{{\ot k}}. \tag 1.7 $$

\m The centralizer algebra $\End_{\g}(M)$ of transformations commuting
with the action of $\g = \fgl_\nn$ on $M = V^{\ot k} \ot (V^*)^{\ot k} $
was investigated in [BCHLLS], where it was shown to be a homomorphic
image of a certain algebra $B_{k,k}(\nn)$ of diagrams with walls.  A
diagram in $B_{k,k}(\nn)$ consists of two rows of vertices with $2k$
vertices in each row.  There is a wall separating the first $k$
vertices on the left in each row from the $k$ vertices on the right.
Each vertex is connected to precisely one edge but with the
requirement that horizontal edges must cross the wall, but vertical
edges cannot cross.  The product $d_1 d_2$ of two diagrams $d_1$ and
$d_2$ is obtained by placing $d_1$ above $d_2$, identifying the bottom
row of $d_1$ with the top row of $d_2$, and following the resulting
paths.  Cycles in the middle are deleted, but there is a scalar
factor, which is $\nn$ to the number of middle cycles.  For example, in
$B_{5,5}(\nn)$ we would have the following product,

$$
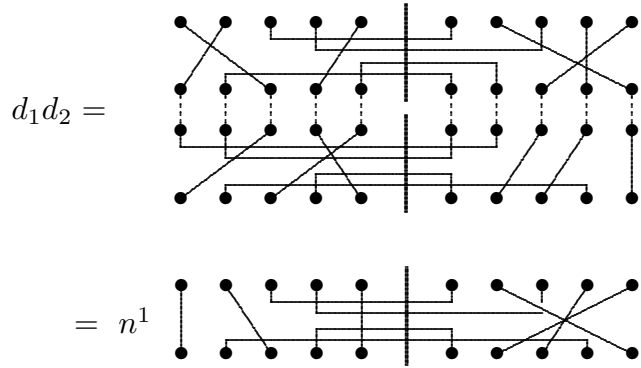
{\beginpicture
 \setcoordinatesystem units <0.5cm,0.3cm>         % sets scale
 \setplotarea x from 0 to 7, y from 0 to 14    % sets plot size up
 \linethickness=0.5pt                          % sets line thickness
 \put{$d_{1}d_{2} =$} at -2.7 6.25
 \put{$\bullet$} at 0 7 \put{$\bullet$} at 0 9.5
 \put{$\bullet$} at 1 7 \put{$\bullet$} at 1 9.5
 \put{$\bullet$} at 2 7 \put{$\bullet$} at 2 9.5
 \put{$\bullet$} at 3 7 \put{$\bullet$} at 3 9.5 
 \put{$\bullet$} at 4 7 \put{$\bullet$} at 4 9.5
 \put{$\bullet$} at 6 7 \put{$\bullet$} at 6 9.5 
 \put{$\bullet$} at 7 7 \put{$\bullet$} at 7 9.5 
 \put{$\bullet$} at 8 7 \put{$\bullet$} at 8 9.5 
 \put{$\bullet$} at 9 7 \put{$\bullet$} at 9 9.5 
 \put{$\bullet$} at 10 7 \put{$\bullet$} at 10 9.5  
 \plot 0 9.5  2 7 /
 \plot 1 9.5  0 7 /
 \plot 2 9.5  2 8.9 /
 \plot 2 8.9  6 8.9 /
 \plot 6 8.9  6 9.5 /
 \plot 3 9.5  3 8.5  /
 \plot 3 8.5   8 8.5 /
 \plot 1  7    1 7.6 /
 \plot 1  7.6   6 7.6 /
 \plot 6 7.6   6  7 /
 \plot 4 7     4 8 /
 \plot 4 8   7 8 /
 \plot 7 8   7 7 /
 \plot 8 8.5   8 9.5  /
 \plot 4 9.5  3 7 /
 \plot 4.98 6.6  4.98 10.2 /
 \plot 5.02 6.6 5.02 10.2 / 
 \plot 7 9.5  10 7 /
 \plot 9 9.5  9  7 /
 \plot 10 9.5  8 7 /
 \put{$\bullet$} at 0 3 \put{$\bullet$} at 0 5.5
 \put{$\bullet$} at 1 3 \put{$\bullet$} at 1 5.5
 \put{$\bullet$} at 2 3 \put{$\bullet$} at 2 5.5
 \put{$\bullet$} at 3 3 \put{$\bullet$} at 3 5.5 
 \put{$\bullet$} at 4 3 \put{$\bullet$} at 4 5.5
 \put{$\bullet$} at 6 3 \put{$\bullet$} at 6 5.5 
 \put{$\bullet$} at 7 3 \put{$\bullet$} at 7 5.5 
 \put{$\bullet$} at 8 3 \put{$\bullet$} at 8 5.5 
 \put{$\bullet$} at 9 3 \put{$\bullet$} at 9 5.5
 \put{$\bullet$} at 10 3 \put{$\bullet$} at 10 5.5 
 \plot 0 5.5  0 5.6    /
 \plot 0 5.8 0  5.9 /
 \plot 0 6.1 0  6.2 /
 \plot 0 6.4 0 6.5 /
 \plot 0 6.7 0 6.8 /
\plot 1 5.5  1 5.6    /
 \plot 1 5.8 1  5.9 /
 \plot 1 6.1 1  6.2 /
 \plot 1 6.4 1 6.5 /
 \plot 1 6.7 1 6.8 /
 \plot 2 5.5  2 5.6    /
 \plot 2 5.8 2  5.9 /
 \plot 2 6.1 2  6.2 /
 \plot 2 6.4 2 6.5 /
 \plot 2 6.7 2 6.8 /
 \plot 3 5.5  3 5.6    /
 \plot 3 5.8 3  5.9 /
 \plot 3 6.1 3  6.2 /
 \plot 3 6.4 3 6.5 /
 \plot 3 6.7 3 6.8 /
 \plot 4 5.5  4 5.6    /
 \plot 4 5.8 4  5.9 /
 \plot 4 6.1 4  6.2 /
 \plot 4 6.4 4 6.5 /
 \plot 4 6.7 4 6.8 /
  \plot 4.98 2.5  4.98 6.1 /
 \plot 5.02 2.5 5.02 6.1 / 
 \plot 6 5.8 6 5.9 /
 \plot 6 6.1 6  6.2 /
 \plot 6 6.4 6 6.5 /
 \plot 6 6.7 6 6.8 /
 \plot 7 5.8 7 5.9 /
 \plot 7 6.1 7  6.2 /
 \plot 7 6.4 7 6.5 /
 \plot 7 6.7 7 6.8 /
 \plot 8 5.8 8 5.9 /
 \plot 8 6.1 8  6.2 /
 \plot 8 6.4 8 6.5 /
 \plot 8 6.7 8 6.8 /
 \plot 9 5.8 9 5.9 /
 \plot 9 6.1 9  6.2 /
 \plot 9 6.4 9 6.5 /
 \plot 9 6.7 9 6.8 /
 \plot 10 5.8 10 5.9 /
 \plot 10 6.1 10 6.2 /
 \plot 10 6.4 10 6.5 /
 \plot 10 6.7 10 6.8 / 
 \plot 0  5.5  0 4.9 /
 \plot 0  4.9  7 4.9 /
 \plot 7  4.9  7 5.5 /
 \plot 1  5.5  1 4.5 /
 \plot 1  4.5  6 4.5 /
 \plot 6 4.5   6 5.5 /
 \plot 2 5.5   0 3 /
 \plot 3 5.5   4 3 /
 \plot 4 5.5   2 3 /
 \plot 1 3     1 3.5 /
 \plot 1 3.5   9 3.5 /
 \plot 9 3.5   9 3 /
 \plot 3 3.9   3 3 /
 \plot 3 3.9   6 3.9 /
 \plot 6 3.9  6 3 /
 \plot 10 5.5 10 3 /
 \plot 9 5.5  8 3 /
 \plot 8 5.5  7 3 /
 \endpicture}$$
 \vskip -.9 truein 
 %%%%%%%%%%%%%%%%%%%%%%%%%%%%%%%%%%%%%%%%%%%%%%%%%%%
 $${\beginpicture
 \setcoordinatesystem units <0.5cm,0.3cm>         % sets scale
 \setplotarea x from 0 to 7, y from 0 to 14    % sets plot size up
 \linethickness=0.5pt                          % sets line thickness
 \put{$ \qquad = \ \nn^1$ \qquad} at -1.4 8.25
 \put{$\bullet$} at 0 7 \put{$\bullet$} at 0 9.5
 \put{$\bullet$} at 1 7 \put{$\bullet$} at 1 9.5
 \put{$\bullet$} at 2 7 \put{$\bullet$} at 2 9.5
 \put{$\bullet$} at 3 7 \put{$\bullet$} at 3 9.5 
 \put{$\bullet$} at 4 7 \put{$\bullet$} at 4 9.5
 \put{$\bullet$} at 6 7 \put{$\bullet$} at 6 9.5 
 \put{$\bullet$} at 7 7 \put{$\bullet$} at 7 9.5 
 \put{$\bullet$} at 8 7 \put{$\bullet$} at 8 9.5 
 \put{$\bullet$} at 9 7 \put{$\bullet$} at 9 9.5 
 \put{$\bullet$} at 10 7 \put{$\bullet$} at 10 9.5  
 \plot 0 9.5  0 7 /
 \plot 1 9.5  2 7 /
 \plot 2 9.5  2 8.9 /
 \plot 2 8.9  6 8.9 /
 \plot 6 8.9  6 9.5 /
 \plot 3 9.5  3 8.5  /
 \plot 3 8.5   8 8.5 /
 \plot 1 7     1 7.5 /
 \plot 1 7.5   9 7.5 /
 \plot 9 7.5   9 7 /
 \plot 3 7.9   3 7 /
 \plot 3 7.9   6 7.9 /
 \plot 6 7.9  6 7 / 
 \plot 8 9.5   8 8.9  /
 \plot 4 9.5  4 7 /
 \plot 4.98 6.6  4.98 10.2 /
 \plot 5.02 6.6 5.02 10.2 / 
 \plot 7 9.5  10 7 /
 \plot 9 9.5  8  7 /
 \plot 10 9.5  7 7 /
 \endpicture}$$

\vskip -.5 truein The group $S_k \times S_k$ acts on $M$, where the
first copy of the symmetric group $S_k$ acts on the first $k$ factors
and the second copy on the next $k$ factors by place permutation.
These actions commute with the $\g$-action, and so afford
transformations in $\End_{\g}(M)$.  There is a representation $\phi\,:
\, B_{k,k}(\nn) \rightarrow \End_{\g}(M)$ of the algebra
$B_{k,k}(\nn)$ on $M$ which commutes with the $\g$-action.  Under this
representation, the diagrams in $B_{k,k}(\nn)$ having no horizontal
edges are mapped to the place permutations coming from $S_k \times
S_k$.  The identity element in $B_{k,k}(\nn)$ is just the diagram with
each node in the top row connected to the one directly below it in the
second row, and it maps to the identity transformation in
$\End_{\g}(M)$.  Under $\phi$, a diagram such as the one pictured
below is mapped to a contraction mapping (in this case to $c_{3,1}$).
\smallskip

$$\beginpicture
 \setcoordinatesystem units <0.5cm,0.3cm>         % sets scale
 \setplotarea x from 0 to 10, y from -1.25 to 1.25   % sets plot size up
 \linethickness=0.5pt                          % sets line  
 \put{$\bullet$} at 0 -1.25 \put{$\bullet$} at 0 1.25
 \put{$\bullet$} at 1 -1.25 \put{$\bullet$} at 1 1.25
 \put{$\bullet$} at 2 -1.25 \put{$\bullet$} at 2 1.25
 \put{$\bullet$} at 3 -1.25 \put{$\bullet$} at 3 1.25 
 \put{$\bullet$} at 4 -1.25 \put{$\bullet$} at 4 1.25
 \put{$\bullet$} at 6 -1.25 \put{$\bullet$} at 6 1.25 
 \put{$\bullet$} at 7 -1.25 \put{$\bullet$} at 7 1.25 
 \put{$\bullet$} at 8 -1.25 \put{$\bullet$} at 8 1.25
 \put{$\bullet$} at 9 -1.25 \put{$\bullet$} at 9 1.25 
 \put{$\bullet$} at 10 -1.25 \put{$\bullet$} at 10 1.25 
 \plot 0 1.25   0 -1.25 /
 \plot 1 1.25   1 -1.25 /
 \plot 2 1.25   2 .65 /
 \plot 2 .65  6 .65  /
 \plot 6 .65  6 1.25 /
 \plot 4 1.25   4 -1.25 /
 \plot 7 1.25    7 -1.25 /
 \plot 3 1.25   3 -1.25 /
 \plot 8 1.25    8 -1.25 / 
 \plot 9 1.25    9 -1.25 / 
 \plot 10 1.25   10 -1.25 /
 \plot 2 -.65  2 -1.25 /
 \plot 2 -.65  6 -.65 /
 \plot 6 -.65 6 -1.25 /  
 \plot 4.98 -1.95  4.98 1.95 /
 \plot 5.02 -1.95 5.02 1.95 / 
 \endpicture  \tag 1.8$$ 
  
\smallskip

  It is shown in [BCHLLS] that the 
algebra $\End_{\g}(M)$  is generated  by $S_k \times S_k$ and the
contraction maps $c_{i,j}$, and the above mapping $\phi$ is
an isomorphism if $\nn \geq 2k$.     Moreover  [BCHLLS]  describes
the projection maps onto the irreducible summands of $M$ in the following
way.   
\m
Suppose for some integer $\rr$ satisfying $0 \leq \rr \leq k$  that 
$\u s = \{s_1,\dots,s_{k-\rr}\}$  and $\u t = \{t_1,\dots, t_{k-\rr}\}$
are ordered subsets of $\{1,\dots,k\}$ of cardinality $k-\rr$ with
$s_1< s_2 < \dots < s_{k-r}$, 
and define the following product
$$c_{\u s,\u t} \eqdef c_{s_1,t_1} \cdots c_{s_{k-\rr},t_{k-\rr}} \tag 1.9$$
of the  contraction maps $c_{s_{i},t_{i}}$.  Then $c_{\u s,\u t}$
belongs to the centralizer algebra $\End_{\g}(M)$.  There is
a corresponding product of diagrams in $B_{k,k}(\nn)$ like the one
displayed in (1.8), which $\phi$ maps onto $c_{\u s,\u t}$.
\m
Assume $\lambda = (\lambda_1 \geq \lambda_2 \geq \cdots)$ is a 
partition of $\rr$.  Associated to $\lambda$
is its Young frame or Ferrers diagram having $\lambda_i$ boxes in
the $i$th row.  A standard tableau is a filling of the boxes in the
diagram of $\lambda$ in such a way that the entries increase 
from left to right across each row and down each column. 
Let $T$ be a standard tableau of shape $\lambda$ with
entries in $\u s^{c} = \{1,\dots,k\}\setminus \{s_1,\dots,s_{k-\rr}\}$.  
Associated to $T$ is its Young symmetrizer 

$$y_T = \Big( \sum_{\rho \in R_T} \rho \Big)\Big( \sum_{\gamma \in C_T} 
\sgn(\gamma) \gamma \Big),  \tag 1.10$$

\n where the first sum ranges over the row group $R_T$  of $T$, which
consists of all permutations in $S_k$ that transform each entry of $T$
to an entry in the same row, and the second sum is over the column
group $C_T$  of $T$ of permutations that move each entry of $T$ to an entry
in the same column.  For example,

$$y_{{\beginpicture
\setcoordinatesystem units <0.3cm,0.3cm>         % sets scale
\setplotarea x from 0 to 1, y from 0 to 2    % sets plot size up
\linethickness=0.3pt              % sets line thickness
\putrule from 0 2 to 2 2        %  draws horizontal lines         
\putrule from 0 1 to 2 1         %           
\putrule from 0 0 to 1 0          %
\putrule from 0 0 to 0 2         %
\putrule from 1 0 to 1 2          %
\putrule from 2 1 to 2 2  %
\put 1 at .5 1.5 % 
\put 5 at 1.5 1.5 %
\put 4 at .5 .5 %
\endpicture}  } = (\id + (1\, 5))(\id - (1\, 4)),$$

\n which belongs to the group algebra $\Bbb C S_k$ of the symmetric
group $S_k$.  The map $y_T$ is an {\it essential idempotent}, that is,
there is an integer $m$ so that $y_T^2 = m y_T$.

\m Similarly, assume for some partition $\mu \vdash \rr$ that $T^{*}$ is
a standard tableau of shape $\mu$ with entries chosen from $\u t^{c} =
\{1,\dots,k\} \setminus \{t_1,\dots, t_{k-\rr}\}$.  The mapping

$$ y_Ty_{T^*}c_{\u s,\u t} \tag 1.11$$

\n is an essential idempotent in $\End_{\g}(M)$. (Note that here we
are supposing that $y_T$ acts on the factors in $V^{\ot k}$ and
$y_{T^{*}}$ on the factors in $(V^{*})^{\ot k}$ by place permutations,
and that $\id$ is the identity map on $V^{\ot k}$ or $(V^{*})^{\ot
k}$, respectively.)  Moreover, the
collection of all maps $y_Ty_{T^*}c_{\u s,\u t} \bigl($ as $\rr = 0,1,\dots,k$;
\ $\u s,\u t$ range over all possible choices of ordered subsets of
cardinality $k-\rr$ in $\{1,\dots,k\}$; \ $\lambda$ and $\mu$ range over all
partitions of $\rr$; and $T$ (resp. $T^*$) ranges over all standard
tableaux of shape $\lambda$ (resp. $\mu$) with entries in $\u s^{c}$
(resp. in $\u t^{c}) \bigr)$  gives all the projections onto the irreducible
summands of $M$ (this can be found in [BCHLLS]).

\medskip 
Now for the idempotent $e$  in (1.6) we may apply the standard result,

$$\End_{\g}(\fsl_\nn^{\ot k}) \cong \End_{\g}(eM) =
e \End_{\g}(M) e \mid_{eM}, \tag 1.12$$
(see for example, [CR, Lemma 26.7] or [BBL, Prop. 1.1]).  
\b
\proclaim {Lemma 1.13}  Assume
$y = y_Ty_{T^*}c_{\u s,\u t}$.
If $c_{\u s,\u t}$  contains one of the contraction maps
$c_{j,j}$ for some $j=1,\dots, k$, then $ey= 0 = ye$.
\endproclaim

\b
\demo{Proof}  The mappings $y_T, y_{T^*}, c_{s_i,t_i}, \, 
i=1,\dots,k-\rr$,
all commute with one another as they operate on different tensor factors.
If one of the contraction maps in $y$ equals $c_{j,j} = \nn p_j$,
then moving it to the far right produces a product $p_je =
p_j (\id-p_j) \prod_{\ell \neq j} (\id-p_{\ell}) = 0$ in $ye$,
so $ye = 0$.   The argument for $ey$ is similar. 
\qed  
\enddemo
\b 
In [BCHLLS, Def. 2.4] (compare also [H1]) a certain simple tensor 
$x_{T,T^{*},\u s, \u t} =
u_1 \ot \cdots \ot u_k \ot w_1^{*}\ot \cdots w_{k}^{*}$  of $M$ is constructed 
via the algorithm 
$$\aligned 
u_{p} & = \cases v_1 & \hskip .5 truein  \text{if}\ p \in \u s \\
v_j & \hskip .5 truein  \text{if  $p \in \u {s}^{c}$ \ and $p$ is in the $j$th row of $T$} 
\endcases \\
w_{p}^* & = \cases v_1^* & \quad \text{if} \ p \in \u t \\
v_{\nn-j+1}^{*} & \quad \text{if  $p \in \u {t}^{c}$ \ and $p$ is in the $j$th 
row of $T^*$} 
\endcases \endaligned \tag 1.14$$

\n When $y = y_Ty_{T^*}c_{\u s,\u t}$ is applied to the simple tensor $x = x_{T,T^{*},\u s, \u t}$ 
the result $y x$ is a nonzero highest 
weight vector  in $y M$.  Moreover, all the highest weight vectors 
in $M$ are produced in this fashion.
\m
Observe that the factors in $x$ lie in $\{v_1,\dots,v_\rr, v_1^{*},v_\nn^{*}, 
\dots, v_{\nn+1-\rr}^{*}\}$.   When the pair $(s_i,t_i)$ belongs to $(\u s, \u 
t)$, then the vector $v_1$ lies in slot $s_i$ in $V^{\ot k}$, and $v_1^*$ lies in 
slot $t_i$ in $(V^{*})^{\ot k}$.   Replace $v_1$ by $v_{\rr+i}$ and 
$v_1^*$ by $v_{\rr+i}^{*}$ in slots $s_i$ and $t_i$ for $i=1,\dots,k-\rr$,
to produce a new simple tensor $x'$.    Then $y x = y x'$, as the 
effect of applying a contraction  to   $v_1 \ot v_1^{*}$ or to
$v_{\rr+i} \ot v_{\rr+i}^{*}$ is the same.     However, if  
$s_i \neq t_i$ for any $i=1,\dots,k-\rr$, then $p_j x' = 0$ for all 
$j$.  The reason for this is that the vector factors of $x'$  
belong to  $\{v_1,\dots,v_\rr,v_{\rr+1},\dots,v_k,v_\nn^{*}, 
\dots, v_{\nn+1-\rr}^{*}, v_k^{*}, \dots, v_{\rr+1}^{*}\}$.  If
$\nn \geq 2k$, these are all distinct.  As $s_i \neq t_i$ for any 
$i=1,\dots,k-\rr$,
slot $j$ on the left and slot $j$ on the right do not contain a pair
of dual vectors (of the form $v_\ell, v_{\ell}^{*}$).    Therefore $p_j x' = 0$ for all $j$ and $e x' = x'$.   
\m
In  [BCHLLS, Thm.~2.5]  it is shown that 
$y_Ty_{T^*}c_{\u s,\u t}\,x_{T,T^{*},\u s, \u t}$  
is a maximal vector in $M$ of highest weight 
$(\lambda,\mu)$, where $\lambda$ is
the shape of $T$ and $\mu$ is the shape of $T^*$.  The $\g$-module
$U(\g)y_Ty_{T^*}c_{\u s,\u t}\,x_{T,T^{*},\u s, \u t}$ generated
by that vector (where $U(\g)$
is the universal enveloping algebra of $\g$)
is isomorphic to the irreducible $\g$-module $L(\lambda,\mu)$
with highest weight $(\lambda,\mu)$.    Moreover, by [BCHLLS, Thm.~2.11], 
the decomposition of $M$ into 
irreducible $\g$-modules is given by 

$$M = \bigoplus U(\g) \, y_Ty_{T^*}c_{\u s,\u t}\, x_{T,T^{*},\u s, \u t}, \tag 1.15$$
where the sum is over all $T,{T^*},\u s,\u t$  as  $\rr = 0,1,\dots,k$;
\ $\u s,\u t$ range over all possible choices of ordered subsets of
cardinality $k-\rr$ in $\{1,\dots,k\}$; \ $\lambda$ and $\mu$ range over all
partitions of $\rr$; and $T$ (resp. $T^*$) ranges over all standard
tableaux of shape $\lambda$ (resp. $\mu$) with entries in $\u s^{c}$
(resp. in $\u t^{c})$.   Since $ey_Ty_{T^*}c_{\u s,\u t} = 0$
whenever $c_{\u s,\u t}$ contains a pair $c_{s_i,t_i}$ with $s_i=t_i$  by Lemma 1.13, 
and since $y_Ty_{T^*}c_{\u s,\u t}\, x_{T,T^{*},\u s, \u t}
= y_Ty_{T^*}c_{\u s,\u t}\, x_{T,T^{*},\u s, \u t}' =
y_Ty_{T^*}c_{\u s,\u t}\, e x_{T,T^{*},\u s, \u t}'$, we have
the following:

\b
\proclaim {Proposition 1.16} Assume $\nn \geq 2k$.  
Then 
$$\aligned
eM  &= \sum_{T,T^*,\u s, \u t}  U(\g) \, e \, y_Ty_{T^*}c_{\u s,\u t}\, x_{T,T^{*},\u s, \u t}\\
&= \sum_{T,T^*,\u s, \u t} U(\g) \, e\, y_Ty_{T^*}c_{\u s,\u t} \, e\,  x_{T,T^{*},\u s, \u t}',  \endaligned $$
where $s_i \neq t_i$ for any  pair $(s_i,t_i)$ in $(\u s,\u t)$.
 \endproclaim 
 
\m
Assume $y = y_Ty_{T^*}c_{\u s,\u t}$ is such that $s_i \neq t_i$ for 
any $i=1,\dots,k-\rr$, and let $x' = x_{T,T^*,\u s,\u t}'$ be the vector
constructed above.    Consider the $U(\g)$-module map
$U(\g)yx' \emap e U(\g)y x' = U(\g)eyx'$.      
Since $U(\g)yx'$ is an irreducible $\g$-submodule of $M$, this map is 0
or an isomorphism.       Now 

$$eyx' = \sum_{J \subseteq \{1,\dots,k\}} (-1)^{|J|} p_J y x' = yx' + \sum_{J \neq \emptyset}  (-1)^{|J|} p_J y x',$$

\noindent where $p_J = \prod_{j \in J} p_j$ as before.  
The right-hand sum is a linear combination of simple tensors
$v_{\ell_1}\ot \cdots \ot v_{\ell_k} \ot v_{m_1}^* \ot \cdots \ot v_{m_k}^*$.
The simple tensor $x'$ does not occur among them, because the
map $p_j$ (for $j=1,\dots,k$) places  $v_i$ in slot $j$ on the
left and $v_i^*$ in slot $j$ on the right, and $x'$ has no such dual pairs
in those particular slots for any $j=1,\dots,k$.    But $x'$ occurs in $yx'$ with coefficient
equal to $|R_T| |R_{T^*}|$, the product of orders of the row groups
of $T$ and $T^*$.   Consequently, since the simple tensors form a basis for $M$, 
we have $eyx' \neq 0$.   Thus,  the above map $E$ is an isomorphism, and $eU(\g)yx' = U(\g)eyx'$
is an irreducible $\g$-module isomorphic to $L(\lambda,\mu)$.  
We have
proved part (1) of the following:
\b
\proclaim {Theorem 1.17}  Assume $\nn \geq 2k$, $\g = \fgl_\nn$, and $M
= V^{\ot k} \ot (V^*)^{\ot k}$. 
\roster
\item  Let $y = y_Ty_{T^*}c_{\u s,\u t}$, where $s_i \neq t_i$ for
any $i=1,\dots,k-\rr$,  and assume $x' = x_{T,T^*,\u s, \u t}' =
u_1' \ot \cdots \ot u_k' \ot (w_1^*)' \ot \cdots \ot  (w_k^*)'$ where 
 
 $$\aligned 
u_{p}' & = \cases v_{r+i}  & \hskip .45 truein  \text{if}\ p = s_i \\
v_j & \hskip .45 truein  \text{if  $p \in \u {s}^{c}$ \ and $p$ is in the $j$th row of $T$} 
\endcases \\
(w_{p}^*)' & = \cases v_{r+i}^* & \ \quad \text{if}  \ p = t_i \\
v_{\nn-j+1}^{*} & \quad \ \,  \text{if  $p \in \u {t}^{c}$ \ and $p$ is in the $j$th 
row of $T^*$.} 
\endcases \endaligned$$
 
\noindent Then $eye U(\g)x' = 
U(\g)eyx' = e U(\g)yx'$ is an irreducible $\g$-submodule of $eM$ of highest weight
$(\lambda,\mu)$ where $\lambda$ is the shape of $T$ and
$\mu$ is the shape of $T^*$.
\item $\fsl_\nn^{\ot k} \cong e M = \bigoplus  U(\g)eyx'$, where the 
sum is over all  $T,{T^*},\u s,\u t$ as $\rr = 0,1,\dots,k$;
\ $\u s,\u t$ range over all possible choices of ordered subsets of
cardinality $k-\rr$ in $\{1,\dots,k\}$ such that $s_i \neq t_i$ for any $i$;  \ $\lambda$ and $\mu$ range over all
partitions of $\rr$; and $T$ (resp. $T^*$) ranges over all standard
tableaux of shape $\lambda$ (resp. $\mu$) with entries in $\u s^{c}$ (resp. in $\u t^{c}$).  
\endroster  \endproclaim 

\b
\demo {Proof}   From Proposition 1.16 we know that

$$eM = \sum U(\g) eyx',  \tag 1.18$$

\noindent where the sum is over all $yx'$  with $y = y_Ty_{T^*}c_{\u s,\u t}$
 and $x' = x_{T,T^*,\u s, \u t}'$, 
as \ $\u s,\u t$ range over all possible choices of ordered subsets of
cardinality $k-\rr$ in $\{1,\dots,k\}$ for $\rr = 0,1,\dots,k$ such that $s_i \neq t_i$ for any $i$;  \ $\lambda$ and $\mu$ range over all
partitions of $\rr$; and $T$ (resp. $T^*$) ranges over all standard
tableaux of shape $\lambda$ (resp. $\mu$) with entries in $\u s^{c}$
(resp. in $\u t^{c}$).     What remains to be
 shown is the sum is direct.
 We have argued previously that  the map, 
 
 $$E: U(\g)yx' \rightarrow eU(\g)yx' = U(\g)eyx'   $$
 
\noindent given by restricting $e$  to $U(\g)yx'$
is an isomorphism of
$\g$-modules for $y = y_Ty_{T^*}c_{\u s,\u t}$ and $x'  = x_{T,T^*,\u s, \u t}'$, 
such that $s_i \neq t_i$ for any $i$.
Fix one such idempotent $y_{\star}$ and consider the intersection 

$$U(\g)e y_{\star} x_\star' \, \, \textstyle{\bigcap} \, \, \displaystyle{\sum_{y \neq y_\star } U(\g)eyx' }
= eU(\g)y_\star x_\star' \, \, \textstyle{\bigcap} \, \, \displaystyle{ e \Bigg(\sum_{y \neq y_\star} U(\g)yx'\Bigg)}$$

\n of $U(\g)e y_\star x_\star'$ with the sum over the remaining ones.   Then

$$eU(\g)y_\star x_\star'\, \, \textstyle{\bigcap} \, \,   \displaystyle{e\Bigg(\sum_{y \neq y_\star} 
U(\g)yx'\Bigg)} \ \einv  U(\g)y_\star x_\star' \, \, \textstyle{\bigcap} \, \, \displaystyle{\sum_{y \neq y_\star } U(\g)yx'} .$$

\n But $U(\g)y_\star x_\star' \, \, \bigcap \, \,  \sum_{y \neq y_\star} U(\g)yx'  = 0$  by (1.15). 
Thus, the sum in (1.18) is direct and we have (2).  \qed \enddemo

\m
\b
\head \S 2.  Multiplicities \endhead 
\b
Knowing that 

$$\fsl_\nn^{\ot k} \cong e M = \bigoplus \, U(\g)\,eyx', $$

\n 
where the sum is
over all $yx' = y_Ty_{T^*}c_{\u s,\u t}\,x_{T,T^*,\u s, \u t}'$ 
such that $s_i \neq t_i$ for any $i$, we may deduce 
 the multiplicity of a particular
irreducible summand in $\fsl_\nn^{\ot k}$ labelled
by $(\lambda,\mu)$, where $\lambda,\mu \vdash \rr$ and $\rr = 0,1,\dots,k$.
That multiplicity is the
number of $(T,T^*,\u s, \u t)$ with $T$ having shape
$\lambda$, \ $T^*$ having shape
$\mu$,  and $c_{\u s,\u t}$ having no  pairs $s_i = t_i$.  
 \m
 
Counting the number of $c_{\u s, \u t}$ with at
least $j$ factors of the form $c_{\ell,\ell}$, we have $
\displaystyle{k \choose j}$
for the choice of those contractions, $\displaystyle{k-j \choose k-\rr-j}$
choices for the remaining $s_i$'s in $\u s$, and 
$\displaystyle {k-j \choose k-\rr-j}$ for the rest of the $t_i$'s  
in $\u t$, and $(k-\rr-j)!$ for the number of ways to 
pair the chosen $s_i$'s with the chosen $t_i$'s.  
Thus, the number of such
$y_Ty_{T^*}c_{\u s,\u t}$  with at least $j$ contractions
of the form $c_{\ell,\ell}$
is 

$${k \choose j}{k-j \choose k-\rr-j}^2 (k-\rr-j)! f^\lambda f^\mu
= {k \choose j}{k-j \choose \rr}^2 (k-\rr-j)! f^\lambda f^\mu, \tag 2.1$$ 
where $f^\la$ (resp. $f^{\mu}$) is the number of standard tableaux of shape $\la$,
(resp. $\mu$). Therefore, by the inclusion-exclusion principle, we have the following 
result.
\b
\proclaim {Theorem 2.2} When $\nn \geq 2k$,  the
multiplicity $m_{\lambda,\mu}^{k}$ in  $\fsl_\nn^{\ot k}$  of the irreducible
module $L(\lambda,\mu)$ for $\g = \fgl_\nn$  with
highest weight $(\lambda,\mu)$, where $\lambda,\mu \vdash r$, is

$$m_{\lambda,\mu}^{k} = f^\la f^\mu \Bigg(\sum_{j = 0}^{k-\rr}(-1)^j
{k \choose j}{k-j\choose \rr}^2 (k-\rr-j)!\Bigg). \tag 2.3$$
\endproclaim 
\b
For a partition $\la$ of $\rr$,  the number
$f^\la$ of standard tableaux of shape $\lambda$ is
given by the well-known hook length formula  

$$f^\la = \frac{\rr!}{h(\la)},$$

\n where $h(\la) = \prod_{(i,j) \in \la}h_{i,j},$ the product
of the {\it hook lengths} of the boxes of $\lambda$.   Thus, $h_{{i,j}}$
is the number of boxes in the $(i,j)$ hook of $\lambda$:
the number of boxes 
 to the right of $(i,j)$
plus the number of boxes below $(i,j)$ plus 1. 
\m 
As a result,  the expression for the multiplicity of the summand
labelled by $(\la,\mu)$ also can be written as 

$$m_{\la,\mu}^{k} = \frac{1}{h(\la)h(\mu)} \sum_{j = 0}^{k-\rr}
(-1)^j \frac{ k! (k-j)!}{j! (k-\rr-j)!}.  \tag 2.4$$
\m
Let us consider
a few interesting special cases.  The multiplicity of the trivial $\g$-module
 in $\fsl_\nn^{\ot k}$ (that is, the dimension of
 the space of $\g$-invariants) is
 
$$m_{\emptyset,\emptyset}^k =  
\sum_{j = 0}^{k}(-1)^j
{k \choose j} (k-j)! = k!\sum_{j=0}^{k}(-1)^j \frac{1}{j!} = D_k,
\tag 2.5$$

\n which is the number of {\it derangements} on the set $\{1,\dots,k\}$
(permutations with no fixed elements).   For small values of
$k$, this number is given by
\smallskip 

$$\matrix  k  &  1 & 2 & 3 & 4 & 5 & 6 & 7 & 8 \\
D_k & 0 & 1 & 2 & 9 & 44 & 265 & 1854 & 14,833 \endmatrix \tag 2.6$$
\smallskip 

\m Next, we compute the number of times the adjoint module
$\frak{sl}_\nn = L({\beginpicture
\setcoordinatesystem units <0.3cm,0.3cm>         % sets scale
\setplotarea x from 0 to 1, y from 0 to 1    % sets plot size up
\linethickness=0.3pt              % sets line thickness
\putrule from 0 1 to 1 1         %  draws horizontal lines         
\putrule from 0 0 to 1 0          %           
\putrule from 0 0 to 0 1          %
\putrule from 1 0 to 1 1          %
\endpicture},{\beginpicture
\setcoordinatesystem units <0.3cm,0.3cm>         % sets scale
\setplotarea x from 0 to 1, y from 0 to 1    % sets plot size up
\linethickness=0.3pt              % sets line thickness
\putrule from 0 1 to 1 1         %  draws horizontal lines         
\putrule from 0 0 to 1 0          %           
\putrule from 0 0 to 0 1          %
\putrule from 1 0 to 1 1          %
\endpicture})$
 occurs in $\frak 
{sl}_\nn^{\otimes k}$.  Using the   
fact that $\fsl_\nn$ is self-dual as a $\g$-module, we
see that the number of times $\frak {sl}_\nn$ appears in $\frak 
{sl}_\nn^{\otimes k}$ is the number of times the trivial module appears
in $\frak 
{sl}_\nn^{\otimes k} \otimes \frak{sl}_\nn = \frak{sl}_\nn^{\otimes (k+1)}$.
Hence, the number of times $\frak {sl}_\nn$ appears in $\frak 
{sl}_\nn^{\otimes k}$ is

$$m_{{\beginpicture
\setcoordinatesystem units <0.3cm,0.3cm>         % sets scale
\setplotarea x from 0 to 1, y from 0 to 1    % sets plot size up
\linethickness=0.3pt              % sets line thickness
\putrule from 0 1 to 1 1         %  draws horizontal lines         
\putrule from 0 0 to 1 0          %           
\putrule from 0 0 to 0 1          %
\putrule from 1 0 to 1 1          %
\endpicture},{\beginpicture
\setcoordinatesystem units <0.3cm,0.3cm>         % sets scale
\setplotarea x from 0 to 1, y from 0 to 1    % sets plot size up
\linethickness=0.3pt              % sets line thickness
\putrule from 0 1 to 1 1         %  draws horizontal lines         
\putrule from 0 0 to 1 0          %           
\putrule from 0 0 to 0 1          %
\putrule from 1 0 to 1 1          %
\endpicture}}^{k} =  D_{k+1}. \tag 2.7$$
\b
This can also be derived from (2.4) which gives

$$
\aligned
m_{{\beginpicture
\setcoordinatesystem units <0.3cm,0.3cm>         % sets scale
\setplotarea x from 0 to 1, y from 0 to 1    % sets plot size up
\linethickness=0.3pt              % sets line thickness
\putrule from 0 1 to 1 1         %  draws horizontal lines         
\putrule from 0 0 to 1 0          %           
\putrule from 0 0 to 0 1          %
\putrule from 1 0 to 1 1          %
\endpicture},{\beginpicture
\setcoordinatesystem units <0.3cm,0.3cm>         % sets scale
\setplotarea x from 0 to 1, y from 0 to 1    % sets plot size up
\linethickness=0.3pt              % sets line thickness
\putrule from 0 1 to 1 1         %  draws horizontal lines         
\putrule from 0 0 to 1 0          %           
\putrule from 0 0 to 0 1          %
\putrule from 1 0 to 1 1          %
\endpicture}}^{k} &  = \sum_{j=0}^{k-1}(-1)^{j} \frac{k! (k-j)}
{j!} = \sum_{j=0}^{k}(-1)^{j} \frac{k! (k-j)}
{j!} \\
& = k\sum_{j=0}^{k}(-1)^{j} \frac{k!}{j!} +
\sum_{j=1}^{k}(-1)^{j-1} \frac{k! }
{(j-1)!} \\
& = k\sum_{j=0}^{k}(-1)^{j} \frac{k!}{j!} +
k\sum_{j=0}^{k-1}(-1)^{j} \frac{(k-1)! }
{j!} \\
& = k(D_k + D_{k-1}) = D_{k+1}. \endaligned \tag 2.8$$
\smallskip 
\n The last equality in (2.8) is a linear recurrence relation 
satisfied by the derangement numbers (see for example, [B, (6.5)]). 

\b
  For any $\frak {g}$-module $X$,
\b
\n $$X \otimes X^* \cong   \End (X)$$ 

\n where the action on the right is
$(g \cdot \psi)(x) =  g\psi(x) - \psi(gx)$ for
all $g \in \g$, $\psi \in \End (X)$, and $x \in X$.   
Considering the $\g$-invariants on both sides, we see that 

$$(X \otimes X^*)^{\frak g} \cong \End(X)^{\frak g} = 
\End_{\g}(X). \tag 2.9$$
Now applying this to $X = \frak{sl}_\nn^{\otimes k} \cong X^{\ast}$,
we have 
$$\End_{\g}(\frak{sl}_\nn^{\otimes k}) \cong (\frak{sl}_\nn^{\otimes 
2k})^{\frak g} \tag 2.10$$  
\n Consequently,

$$\dim \End_{\g}(\frak{sl}_\nn^{\otimes k}) = 
m_{\emptyset,\emptyset}^{2k} = 
D_{2k}, \tag 2.11$$
  
\n the number of 
derangements
on a set of $2k$ elements.

\b
We conclude by displaying the multiplicities $m_{\lambda,\mu}^{k}$ for $k = 4$.  By double
centralizer theory, it follows that 

$$\dim \End_{\g}(\frak{sl}_\nn^{\otimes k}) = \sum_{\lambda,\mu\,
\vdash \rr \leq k} \big( m_{\lambda, \mu}^{k}\big )^{2}.$$

\n The reader can verify that the squares of the
numbers in the following tables do indeed sum to $D_{8} = 14,833$. 
\vfill \eject

\n {\bf Example}: \quad  $m_{\lambda, \mu}^{4}$:   
\b 
$$\matrix  &{\beginpicture
\setcoordinatesystem units <0.3cm,0.3cm>         % sets scale
\setplotarea x from 0 to 1, y from 0 to 1    % sets plot size up
\linethickness=0.3pt              % sets line thickness
\putrule from 0 4.5 to 4 4.5         %  draws horizontal lines         
\putrule from 0 3.5 to 4 3.5          %           
\putrule from 0 3.5 to 0 4.5         %
\putrule from 1 3.5 to 1 4.5          %
\putrule from 2 3.5 to 2 4.5          %
\putrule from 3 3.5 to 3 4.5          %
\putrule from 4 3.5 to 4 4.5          %
\endpicture}  & {\beginpicture
\setcoordinatesystem units <0.3cm,0.3cm>         % sets scale
\setplotarea x from 0 to 1, y from 0 to 2    % sets plot size up
\linethickness=0.3pt              % sets line thickness
\putrule from 0 4.5 to 3 4.5         %  draws horizontal lines         
\putrule from 0 3.5 to 3 3.5          %           
\putrule from 0 2.5 to 1 2.5          %
\putrule from 0 2.5 to 0 4.5         %
\putrule from 1 2.5 to 1 4.5          %
\putrule from 2 3.5 to 2 4.5          %
\putrule from 3 3.5 to 3 4.5       %
\endpicture} & 
{\beginpicture
\setcoordinatesystem units <0.3cm,0.3cm>         % sets scale
\setplotarea x from 0 to 1, y from 0 to 2    % sets plot size up
\linethickness=0.3pt              % sets line thickness
\putrule from 0 4.5 to 2 4.5         %  draws horizontal lines         
\putrule from 0 3.5 to 2 3.5          %           
\putrule from 0 2.5 to 2 2.5          %
\putrule from 0 2.5 to 0 4.5
\putrule from 1 2.5 to 1 4.5         %
\putrule from 2 2.5 to 2 4.5          % 
\endpicture}
& 
{\beginpicture
\setcoordinatesystem units <0.3cm,0.3cm>         % sets scale
\setplotarea x from 0 to 1, y from 0 to 2    % sets plot size up
\linethickness=0.3pt              % sets line thickness
\putrule from 0 4.5 to 2 4.5         %  draws horizontal lines         
\putrule from 0 3.5 to 2 3.5          %           
\putrule from 0 2.5 to 1 2.5          %
\putrule from 0 1.5 to 1 1.5
\putrule from 0 1.5 to 0 4.5
\putrule from 1 1.5 to 1 4.5         %
\putrule from 2 3.5 to 2 4.5          % 
\endpicture} 
& 
{\beginpicture
\setcoordinatesystem units <0.3cm,0.3cm>         % sets scale
\setplotarea x from 0 to 1, y from 0 to 2    % sets plot size up
\linethickness=0.3pt              % sets line thickness
\putrule from 0 4.5 to 1 4.5         %  draws horizontal lines         
\putrule from 0 3.5 to 1 3.5          %           
\putrule from 0 2.5 to 1 2.5          %
\putrule from 0 1.5 to 1 1.5   %
\putrule from 0 .5 to  1  .5
\putrule from 0 .5 to 0 4.5
\putrule from 1 .5 to 1 4.5         % 
\endpicture}\\
{\beginpicture
\setcoordinatesystem units <0.3cm,0.3cm>         % sets scale
\setplotarea x from 0 to 1, y from 0 to 1    % sets plot size up
\linethickness=0.3pt              % sets line thickness
\putrule from 0 1 to 4 1         %  draws horizontal lines         
\putrule from 0 0 to 4 0          %           
\putrule from 0 0 to 0 1        %
\putrule from 1 0 to 1 1          %
\putrule from 2 0 to 2 1          %
\putrule from 3 0 to 3 1          %
\putrule from 4 0 to 4 1         %
\endpicture} & 1 & 3 & 2 & 3 & 1 \\
 {\beginpicture
\setcoordinatesystem units <0.3cm,0.3cm>         % sets scale
\setplotarea x from 0 to 1, y from 0 to 2    % sets plot size up
\linethickness=0.3pt              % sets line thickness
\putrule from 0 1.5  to 3 1.5        %  draws horizontal lines         
\putrule from 0  .5 to 3  .5          %           
\putrule from 0 -0.5 to 1 -0.5          %
\putrule from 0 -0.5 to 0 1.5         %
\putrule from 1 -0.5 to 1 1.5          %
\putrule from 2  .5 to 2 1.5         %
\putrule from 3  .5 to 3 1.5       %
\endpicture}& 3 & 9 & 6 & 9 & 3 \\
{\beginpicture
\setcoordinatesystem units <0.3cm,0.3cm>         % sets scale
\setplotarea x from 0 to 1, y from 0 to 2    % sets plot size up
\linethickness=0.3pt              % sets line thickness
\putrule from 0 1 to 2 1          %  draws horizontal lines         
\putrule from 0  0 to 2  0          %           
\putrule from 0 -1 to 2 -1          %
\putrule from 0 -1 to 0 1
\putrule from 1 -1 to 1 1         %
\putrule from 2 -1 to 2 1          % 
\endpicture} & 2 & 6 & 4 & 6 & 2 \\
{\beginpicture
\setcoordinatesystem units <0.3cm,0.3cm>         % sets scale
\setplotarea x from 0 to 1, y from 0 to 2    % sets plot size up
\linethickness=0.3pt              % sets line thickness
\putrule from 0 1.5 to 2 1.5         %  draws horizontal lines         
\putrule from 0 .5 to 2  .5          %           
\putrule from 0 -.5 to 1 -.5          %
\putrule from 0 -1.5 to 1 -1.5
\putrule from 0 -1.5 to 0 1.5
\putrule from 1 -1.5 to 1 1.5         %
\putrule from 2  .5 to 2 1.5          % 
\endpicture} & 3 & 9 & 6 & 9 & 6  \\
 {\beginpicture
\setcoordinatesystem units <0.3cm,0.3cm>         % sets scale
\setplotarea x from 0 to 1, y from 0 to 2    % sets plot size up
\linethickness=0.3pt              % sets line thickness
\putrule from -.5 1.5 to .5 1.5        %  draws horizontal lines         
\putrule from -.5 .5 to .5 .5          %           
\putrule from -.5 -.5 to .5 -.5          %
\putrule from -.5 -1.5 to .5 -1.5  %
\putrule from -.5 -2.5 to  .5  -2.5
\putrule from -.5 -2.5 to -.5 1.5
\putrule from  .5 -2.5 to .5 1.5        % 
\endpicture}
& 1 & 3 & 2 & 3 & 1  \endmatrix $$
%%%%%%%%%%%%%%%%%%%%%%%%%%%%%%%%%%%%%%%%%%(end of k = 4)
\m
\b
$$\matrix  &{\beginpicture
\setcoordinatesystem units <0.3cm,0.3cm>         % sets scale
\setplotarea x from 0 to 1, y from 0 to 2    % sets plot size up
\linethickness=0.3pt              % sets line thickness
\putrule from 0 3.5 to 3 3.5         %  draws horizontal lines         
\putrule from 0 2.5 to 3 2.5          %           
\putrule from 0 2.5 to 0 3.5         %
\putrule from 1 2.5 to 1 3.5          %
\putrule from 2 2.5 to 2 3.5          %
\putrule from 3 2.5 to 3 3.5
\endpicture}  & {\beginpicture
\setcoordinatesystem units <0.3cm,0.3cm>         % sets scale
\setplotarea x from 0 to 1, y from 0 to 2    % sets plot size up
\linethickness=0.3pt              % sets line thickness
\putrule from 0 3.5 to 2 3.5         %  draws horizontal lines         
\putrule from 0 2.5 to 2 2.5          %           
\putrule from 0 1.5 to 1 1.5          %
\putrule from 0 1.5 to 0 3.5         %
\putrule from 1 1.5 to 1 3.5          %
\putrule from 2 2.5 to 2 3.5   %
\endpicture} & 
{\beginpicture
\setcoordinatesystem units <0.3cm,0.3cm>         % sets scale
\setplotarea x from 0 to 1, y from 0 to 2    % sets plot size up
\linethickness=0.3pt              % sets line thickness
\putrule from 0 3.5 to 1 3.5         %  draws horizontal lines         
\putrule from 0 2.5 to 1 2.5          %           
\putrule from 0 1.5 to 1 1.5          %
\putrule from 0 0.5 to 1 0.5
\putrule from 0 0.5 to 0 3.5         %
\putrule from 1 0.5 to 1 3.5          % 
\endpicture} \\
{\beginpicture
\setcoordinatesystem units <0.3cm,0.3cm>         % sets scale
\setplotarea x from 0 to 1, y from 0 to 1    % sets plot size up
\linethickness=0.3pt              % sets line thickness
\putrule from  0 1  to 3 1          %  draws horizontal lines         
\putrule from 0  0 to 3 0         %           
\putrule from 0  0 to 0 1       %
\putrule from 1  0 to 1 1          %
\putrule from 2  0 to 2 1          %
\putrule from 3  0 to 3 1 
\endpicture}   & 12 & 24 & 12 \\
{\beginpicture
\setcoordinatesystem units <0.3cm,0.3cm>         % sets scale
\setplotarea x from 0 to 1, y from 0 to 2.75    % sets plot size up
\linethickness=0.3pt              % sets line thickness
\putrule from 0 2 to 2 2        %  draws horizontal lines         
\putrule from 0 1 to 2 1         %           
\putrule from 0 0 to 1 0          %
\putrule from 0 0 to 0 2         %
\putrule from 1 0 to 1 2          %
\putrule from 2 1 to 2 2  %
\endpicture}& 24 & 48  & 24 \\
{\beginpicture
\setcoordinatesystem units <0.3cm,0.3cm>         % sets scale
\setplotarea x from 0 to 1, y from 0 to 2.5    % sets plot size up
\linethickness=0.3pt              % sets line thickness
\putrule from 0 1.5 to 1 1.5         %  draws horizontal lines         
\putrule from 0  .5 to 1  .5          %           
\putrule from 0 -.5 to 1 -.5          %
\putrule from 0 -1.5 to 1 -1.5
\putrule from 0 -1.5 to 0 1.5         %
\putrule from 1 -1.5 to 1 1.5          % 
\endpicture} & 12 & 24 & 12 \endmatrix \qquad \qquad
 \matrix  &{\beginpicture
\setcoordinatesystem units <0.3cm,0.3cm>         % sets scale
\setplotarea x from 0 to 1, y from 0 to 1    % sets plot size up
\linethickness=0.3pt              % sets line thickness
\putrule from 0 2 to 2 2         %  draws horizontal lines         
\putrule from 0 1 to 2 1          %           
\putrule from 0 1 to 0 2          %
\putrule from 1 1 to 1 2          %
\putrule from 2 1 to 2 2     %
\endpicture}  & {\beginpicture
\setcoordinatesystem units <0.3cm,0.3cm>         % sets scale
\setplotarea x from 0 to 1, y from 0 to 2.5    % sets plot size up
\linethickness=0.3pt              % sets line thickness
\putrule from 0 2 to 1 2         %  draws horizontal lines         
\putrule from 0 1 to 1 1          %           
\putrule from 0 0 to 1 0          %
\putrule from 0 0 to 0 2          %
\putrule from 1 0 to 1 2     %
\endpicture} \\
{\beginpicture
\setcoordinatesystem units <0.3cm,0.3cm>         % sets scale
\setplotarea x from 0 to 1, y from 0 to 2    % sets plot size up
\linethickness=0.3pt              % sets line thickness
\putrule from 0 1 to 2 1         %  draws horizontal lines         
\putrule from 0 0 to 2 0          %           
\putrule from 0 0 to 0 1          %
\putrule from 1 0 to 1 1          %
\putrule from 2 0 to 2 1     %
\endpicture}  & 42 & 42 \\
{\beginpicture
\setcoordinatesystem units <0.3cm,0.3cm>         % sets scale
\setplotarea x from 0 to 1, y from -1 to 1.75    % sets plot size up
\linethickness=0.3pt              % sets line thickness
\putrule from 0 1 to 1 1         %  draws horizontal lines         
\putrule from 0 0 to 1 0          %           
\putrule from 0 -1 to 1 -1          %
\putrule from 0 -1 to 0 1          %
\putrule from 1 -1 to 1 1     %
\endpicture} & 42 & 42 \endmatrix \qquad  
\matrix & {\beginpicture
\setcoordinatesystem units <0.3cm,0.3cm>         % sets scale
\setplotarea x from 0 to 1, y from 0 to 1    % sets plot size up
\linethickness=0.3pt              % sets line thickness
\putrule from 0 1.5 to 1 1.5         %  draws horizontal lines         
\putrule from 0 .5 to 1 .5          %           
\putrule from 0 .5 to 0 1.5          %
\putrule from 1 .5 to 1 1.5          %
\endpicture}  \\
{\beginpicture
\setcoordinatesystem units <0.3cm,0.3cm>         % sets scale
\setplotarea x from 0 to 1, y from 0 to 1    % sets plot size up
\linethickness=0.3pt              % sets line thickness
\putrule from 0 1 to 1 1         %  draws horizontal lines         
\putrule from 0 0 to 1 0          %           
\putrule from 0 0 to 0 1          %
\putrule from 1 0 to 1 1          %
\endpicture} & 44\endmatrix \qquad \qquad
\matrix &\emptyset  \\ 
\emptyset & 9 \endmatrix $$
\b 
\b 

\head \S 3. The Centralizer Algebra \endhead

\m Now we consider the centralizer algebra $\Cent =
\End_\g(\fsl_\nn^{\otimes k}) = \End_{\fsl_\nn}(\fsl_\nn^{\otimes k})$ and
its representation theory. As has already been pointed out in (1.12), we have an
isomorphism
$$
\Cent \cong e \End_\g(M) e   \tag 3.1
$$
where $e$ is the idempotent defined in (1.6). We also have a
representation $\phi\,: B_{k,k}(\nn) \to \End(M)$ which commutes with the
$\g$-action on $M$.  Thus the image of this representation lies in the
commuting algebra $\End_\g(M)$.  In [BCHLLS, Thm.~5.8] it was shown
that $\phi$ induces an algebra isomorphism
$$
B_{k,k}(\nn) \cong \End_\g(M)   \tag 3.2
$$
for $\nn\ge 2k$.  
 
 \b
Let $c_j$ denote the diagram in $B_{k,k}(\nn)$ 
corresponding to the contraction $c_{j,j}$, but 
scaled by a factor of $1/\nn$.  Then 
  under  the representation
$\phi\,:\,  B_{k,k}(\nn) \rightarrow \End_{\g}(M)$, 
$c_j$ is sent to $p_j$, and 
$b = \prod_{j=1}^k (1 - c_j)$ is mapped to the idempotent $e$. 
\medskip
Let us consider the subspace $A$ spanned
by the diagrams $d$ having no forbidden pairs.  By a forbidden
pair,  we
mean that the $i$th node on the left 
is connected to the $i$th node on the
right of the wall either in the top or in the bottom row of $d$ 
for some $i=1,\dots,k$.  
\m
We claim that the map $B_{k,k}(\nn) \rightarrow bB_{k,k}(\nn)b$
is injective on the subspace $A$ of diagrams with no forbidden
pairs.  Indeed, $\sum_{d \in A} a_d d \mapsto
\sum_{d \in A} a_d bd b = \sum_{d \in A} a_d d  + f$, where $a_d \in 
\Bbb C$ and 
$f$ is a linear combination of diagrams in $B_{k,k}(\nn)$ having
at least one forbidden pair.  The reason for this is that when 
diagrams are multiplied, the horizontal edges in the top row of the top diagram
and the horizontal edges in the bottom row of the bottom diagram
always appear in the resulting product diagram.     Thus, we obtain the following

\b
\proclaim {Proposition 3.3} Let $\nn\ge 2k$. The map $\phi$ induces an
algebra isomorphism between $b B_{k,k}(\nn) b$ and $\Cent =
\hbox{\rm End}_\g(\fsl_\nn^{\otimes k})$.  Moreover, the set of all elements of
the form $b d b$, as $d$ ranges over all diagrams with no forbidden
pairs, is a basis for $b B_{k,k}(\nn) b$.  \endproclaim

\b
\demo{Proof} The first claim follows from the remarks above, so only
the second claim remains to be proved.  We observe that left (resp.,
right) multiplication by $b$ kills any diagram with a forbidden pair
in its top (resp., bottom) row. Since the diagrams form a basis for
$B_{k,k}(\nn)$, the result follows. \qed \enddemo
\b

The basis statement of Proposition 3.3 provides another
proof of (2.11), that the dimension of the centralizer algebra $\Cent$ is
$D_{2k}$.  Indeed, the diagrams with no forbidden pairs are easily
seen to be in bijective correspondence with the permutations $\sigma$
on the set $\{1,\dots,2k\}$ such that $\sigma(i) \ne i$ for all
$i=1, \dots, 2k$.  This correspondence is given by performing two 
``flips'', which take a walled Brauer diagram to the diagram obtained by 
first interchanging the rightmost $k$ dots in its top and bottom rows 
and then switching corresponding dots on the two sides of the wall on the top 
row while retaining the edges.

\m
Let $\rr\le k$ and let $\lambda,\mu$ be fixed partitions of $\rr$.  In
[BCHLLS] $M_{\lambda,\mu}$ was defined to be the space spanned
by all maximal vectors $yx = yx'$, where
$y=y_Ty_{T^*}c_{\u s,\u t}$,    $x=x_{T,T^*,\u s, \u t}$ (notation of (1.14)), 
and $x' = x_{T,T^*,\u s, \u t}'$ as in Theorem 1.17 \ for all pairs
$\u s = \{s_1, \dots, s_{k-\rr}\}$, $\u t = \{t_1,\dots, t_{k-\rr}\}$ of
ordered subsets of $\{1,\dots,k\}$, and all standard tableaux $T$
(resp., $T^*$) of shape $\lambda$ (resp., $\mu$) with entries from $\u
s^c$ (resp., $\u t^c$).  Moreover, for $\nn\ge 2k$, the
$M_{\lambda,\mu}$ provide a complete set of pairwise nonisomorphic
irreducible modules for the algebra $\End_\g(M)$ (and hence also for
$B_{k,k}(\nn)$). 

\b
\proclaim {Lemma 3.4}  Assume $\nn \geq 2k$ and let $y = y_Ty_{T^*}c_{\u
s,\u t}$, $x'=x_{T,T^*,\u s, \u t}'$ as in Theorem 1.17.  Then $eyx'  \ne
0$ if and only if $s_i \neq t_i$ for all pairs $(s_i,t_i)$ in $(\u s,
\u t)$. Hence $eM_{\lambda,\mu} \ne 0$ precisely when this condition 
can be satisfied, and in that case,  $eM_{\lambda,\mu}$ is the
linear span of all the nonzero $eyx'$, $y$ and $x'$ as above. \endproclaim

\b
\demo {Proof} This follows from results in [BCHLLS], Lemma 1.13, 
and its converse, which 
is in the paragraph before Theorem 1.17.  \qed \enddemo
 
\b It is easy to see that $eM_{\lambda,\mu}=0$ when
$\lambda=\mu=\emptyset$ and $k=1$, for in that case it is impossible to
construct a  $y = y_Ty_{T^*}c_{\u s,\u t}$ satisfying the condition 
$s_i \neq t_i$ for all pairs $(s_i,t_i)$ in $(\u s,
\u t)$.  In all other cases $eM_{\lambda,\mu} \neq 0$ when $n \geq 2k$.

\b
\proclaim {Theorem 3.5} Assume $\nn \ge 2k$.  The collection of all
nonzero $eM_{\lambda,\mu}$ for $\lambda,\mu$ partitions of $\rr$,
$\rr=0,1,\dots,k$, forms a complete set of pairwise nonisomorphic
irreducible modules for the algebra $\Cent\cong b B_{k,k}(\nn)
b$. \endproclaim

\b \demo {Proof} It is well-known that if $u$ is an idempotent in
an algebra $A$, the functor $u(-)$ (sometimes called the Schur
functor; see [G, 6.2]) taking $A$-modules to $uAu$-modules is an
exact covariant functor which maps an irreducible module to either an
irreducible module or zero.  In the particular case that
 $A = B_{k,k}(\nn)$ and $u=b$,
this functor takes the irreducible module $M_{\lambda,\mu}$ to
$bM_{\lambda,\mu} = eM_{\lambda,\mu}$. \qed \enddemo  
\b
\proclaim {Theorem 3.6} Assume $\nn\ge 2k$.  Then as a bimodule 
for $\Cent \times \g$, 
$$
 \fsl_\nn^{\otimes k} \cong eM \cong \bigoplus_{\rr=0}^k
 \bigoplus_{\lambda,\mu \vdash \rr} eM_{\lambda,\mu} \otimes
 L(\lambda,\mu),
$$
where the decomposition is into pairwise nonisomorphic
irreducible modules for $\Cent \times \g$. \endproclaim
\b

\demo {Proof} This follows from the previous results and standard 
double-centralizer theory. \qed \enddemo

\b For $\nn\ge 2k$  the dimension of the irreducible $\Cent$-module
$eM_{\lambda,\mu}$ is given by $m_{\lambda,\mu}^k$ (see Theorem 2.2).

\b 
\b

\b
 
\address 
\newline 
Department of Mathematics, University of Wisconsin, Madison,
Wisconsin 53706   \newline 
benkart@math.wisc.edu
\newline 
\newline  
Department of Mathematical and Computer Sciences, 
Loyola University Chicago, Chicago,
Illinois 60626 \newline
doty@math.luc.edu \endaddress 

\enddocument
\end


\Refs 
\widestnumber\key{BCHLLS} 
\b
\ref \key BBL \by G. Benkart, D.J. Britten, F.W. Lemire,
Projection maps for tensor products of $\frak{gl}(r,\Bbb C)$-representations 
\jour Publ. RIMS, Kyoto \vol 28  \yr 1992 \pages
983--1010  \endref
\m
\ref \key BCHLLS \by G. Benkart, M. Chakrabarti,
T. Halverson, R. Leduc, C. Lee, and J. Stroomer,
Tensor product representations of general linear groups
and their connections with Brauer algebras \jour  J. Algebra 
\vol 166 \yr 1994 \pages 529-567  \endref
\m
\ref \key B  \by R.A. Brualdi  \book Introductory Combinatorics, 3rd ed.   
\publ Prentice Hall \publaddr Englewood Cliffs, N.J. \yr 1999  \endref
\m
\ref \key CR \by C.W. Curtis and I. Reiner  \book Representation Theory of
Finite Groups and Associative Algebras  \vol XI \publ Pure and Applied
Math, Interscience Publ. John Wiley \publaddr New York \yr 1962 
\endref 
\m
\ref \key G  \by J.A. Green \book Polynomial Representations of
$\text{GL}_n$ \bookinfo Lecture Notes in Math. \vol 830
\publ Springer-Verlag \publaddr Heidelberg \yr 1980
\endref
\m
\ref \key H1  \by P. Hanlon, On the construction of the maximal 
vectors in the tensor algebra of $\fgl_n$  
\book Combinatorics and algebra (Boulder, Colo., 1983)
Contemp. Math. \vol 34 \pages 73--80 \publ Amer. Math. Soc.
\publaddr Providence R.I. \yr 1984 \endref
\m
\ref \key H2  \by P. Hanlon, On the decomposition of the tensor algebra
of the classical Lie algebras \jour Adv. in Math. \vol 56 \pages
238--282 \yr 1985 \endref


\endRefs
\enddocument 
\end